\documentclass[11pt,a4paper,reqno]{amsart}

\usepackage{amsfonts}
\usepackage{graphicx}
\usepackage{graphics}
\usepackage{epsfig}
\usepackage{anysize}
\usepackage{lipsum}
\usepackage{wrapfig}
\marginsize{3cm}{3cm}{3cm}{3cm}
\newtheorem{theorem}{Theorem}[section]

\newtheorem{lemma}{Lemma}[section]

\newtheorem{corollary}{Corollary}[section]

\makeatletter

\begin{document}

	\title{An Algorithmic Approach to Antimagic Labeling of Edge Corona Graphs}
	\author{D. Nivedha, S. Devi Yamini}
	\thanks{\noindent Vellore Institute of Technology - Department of Mathematics - Kelambakkam - Vandalur Rd, Rajan Nagar, Chennai, Tamil Nadu - 600127, India\\
		\indent \,\,\, e-mail: nivedha.d2020@vitstudent.ac.in;  https://orcid.org/0000-0002-1292-8629\\
		\indent \,\,\, e-mail: deviyamini.s@vit.ac.in;  https://orcid.org/0000-0002-7866-7230 \\}
	
	\begin{abstract}
		\noindent An antimagic labeling of a graph $G$ is a $1-1$ correspondence between the edge set $E(G)$ and $\lbrace 1,2,...,|E(G)|\rbrace$ in which the sum of the labels of edges incident to the distinct vertices are different. The edge corona of any two graphs $G$ and $H$, (denoted by $G$ $\diamond$ $H$) is obtained by joining one copy of $G$ with $|E(G)|$ copies of H such that the end vertices of $i^{th}$ edge of $G$ is adjacent to every vertex in the $i^{th}$ copy of $H$. In this paper, we provide an algorithm to prove that the following graphs admit an antimagic labeling:\\ 
		$-$ $n$-barbell graph $B_n$, $n\geq3$ \\
		$-$ edge corona of a bistar graph $B_{x,n}$ and a $k$-regular graph $H$ denoted by $B_{x,n}\diamond H$, $x,n\geq 2$\\
		$-$ edge corona of a cycle $C_m$ and $C_n$ denoted by $C_m \diamond C_n$, $m,n\geq3$	
		
		
		
	\end{abstract}
	\maketitle 
	\bigskip
	\bigskip
	%
	\section{Introduction}
	The graphs considered in this paper are simple, finite and undirected. Let $V(G)$ and $E(G)$ denote the vertex set and the edge set of the graph $G$ respectively. For a graph $G=(V,E)$ with $q$ edges, an antimagic labeling is a bijection $f:E(G)\rightarrow \lbrace 1,2,...,q \rbrace$ such that $w(u)=\sum\limits_{e\in E(u)}f(e)$ is distinct for all vertices $u\in V(G)$ where $E(u)$ denotes the set of edges incident on the vertex $u$. A graph is said to be antimagic if it admits an antimagic labeling. Paths, cycles, complete graphs, wheel, stars, complete bipartite graphs\cite{11}, graphs of order $n$ with maximum degree at least $n-3$\cite{12}, toroidal grid graphs\cite{13}, regular graphs\cite{2,1}, trees with some restrictions \cite{9}, caterpillars\cite{16} are few of the graph classes proved to be antimagic in the literature.\\
	\indent
	The antimagic labeling on edge corona graphs has comparatively less results than other graph products like cartesian product, lexicographic product, corona product, etc in the literature. Also, motivated by the two conjectures proposed by  Hartsfield and Ringel which is still open\cite{11}, we focus on antimagicness of $n$-barbell graph and edge corona of few graph classes in this paper.\\
	\textbf{Conjecture 1:}\cite{11}
	Every connected graph other than $K_2$ is antimagic.\\
	\textbf{Conjecture 2:}\cite{11}
	Every tree other than $K_2$ is antimagic.
	
	\section{Preliminaries}	
	Let $G$ and $H$ be two vertex-disjoint graphs. Define $[1,n] = \lbrace 1,2,...,n\rbrace$. Let $d(v)$ denote the degree of a vertex $v$ in the graph $G$. We say that $u\in V(G)$ is complete to a graph $H$ if $u$ is adjacent to all the vertices of $H$. Also, an edge $ab\in E(G)$ is complete to a graph $H$ if the vertices $a$ and $b$ are adjacent to all the vertices of $H$. The edge corona of any two graphs $G$ and $H$, (denoted by $G$ $\diamond$ $H$) is obtained by joining one copy of $G$ with $|E(G)|$ copies of H such that the end vertices of $i^{th}$ edge of $G$ is adjacent to every vertex in the $i^{th}$ copy of $H$. The sum of the labels of all the edges incident to a vertex in a graph $G$ is called the vertex sum (denoted by $w(v)$, $v \in V(G)$). The sum of the labels of some edges (i.e., few edges remain unlabeled) in a graph $G$ is known as the partial vertex sum (denoted by $w'(v)$, $v\in V(G)$). A graph is said to be regular if the degrees of all the vertices are same. The join of two graphs $G$ and $H$ is obtained by making every vertex of $G$ adjacent to all the vertices of $H$. A bistar graph $B_{m,s}$ is obtained by joining the apex vertices (a vertex adjacent to all the vertices of a graph) of two vertex disjoint star graphs $K_{1,m}$ and $K_{1,s}$ for $m\geq1$ and $s\geq1$ respectively\rm\cite{7}. An $n$-barbell graph is obtained by adding an edge between two copies of $K_n$, $n\geq 3$\rm\cite{8}.

	\section{Main results}
\begin{figure}[ht!]
		\centering
		\includegraphics[width=1.2\linewidth, height=0.15\textheight]{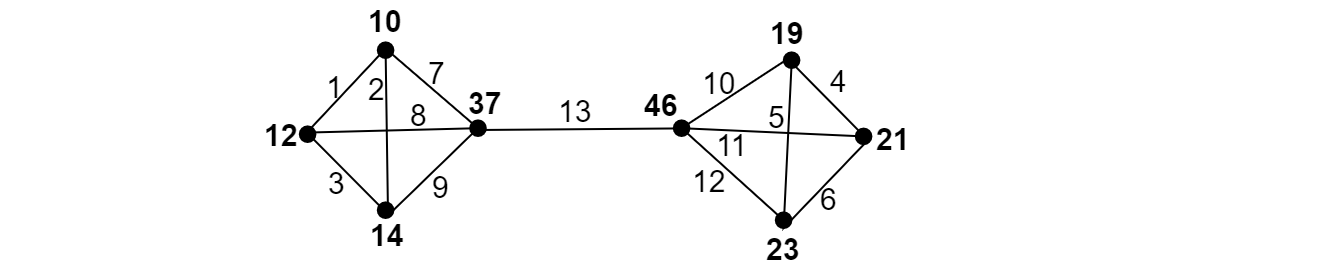}
		\caption{An antimagic labeling of $4$-barbell graph, $B_4$.}
		\label{fig:2}
	\end{figure}

	%
	%
	\begin{theorem}
		The $n$-barbell graph is antimagic for $n\geq3$.\\
		\begin{proof}
			The general representation of an $n$-barbell graph $B_n$ is given in Figure \ref{fig:1}. The dotted line from $u_2$ to $u_{n-2}$ represent a clique $R$ on vertices $\lbrace u_3,u_4,...,u_{n-3}\rbrace$ and the vertices $\lbrace u_1,u_2,u_{n-2},u_{n-1},u_n\rbrace$ are complete to $R$. Similarly, the dotted line from $v_2$ to $v_{n-2}$ represent a clique $S$ on vertices $\lbrace v_3,v_4,...,v_{n-3}\rbrace$ and the vertices $\lbrace v_1,v_2,v_{n-2},v_{n-1},v_n\rbrace$ are complete to $S$. Note that the vertices $\lbrace u_1,u_2,...,u_{n-1}\rbrace$ and $\lbrace v_1,v_2,...,v_{n-1}\rbrace$ need not be in cyclic order. Let $V(B_n)=\lbrace u_1,u_2,...,u_n\rbrace\cup\lbrace v_1,v_2,...,v_n\rbrace$ and $E(B_n)=\lbrace u_iu_j\hspace{0.1cm}\big|\hspace{0.1cm}1\leq i,j\leq n, i\neq j\rbrace\cup\lbrace v_iv_j\hspace{0.1cm}\big|\hspace{0.1cm}1\leq i,j\leq n, i\neq j\rbrace\cup\lbrace uv\rbrace$.\\		
			\begin{figure}[ht!]
				\centering
				\includegraphics[width=1\linewidth, height=0.2\textheight]{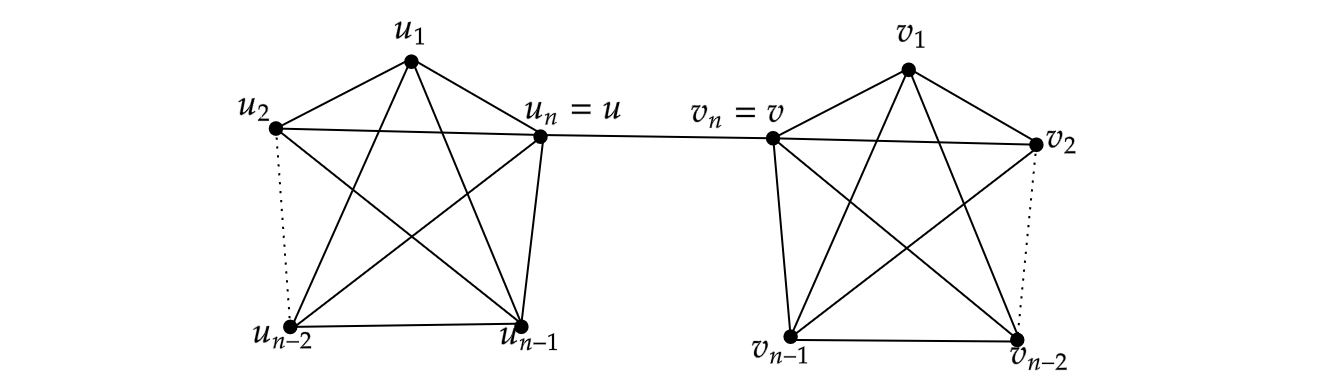}
				\caption{A representation of an $n$-barbell graph $B_n$.}
				\label{fig:1}
			\end{figure}
			\noindent
			The graph $B_n$ contains $2n$ vertices and $n(n-1)+1$ edges. The edge labels are $\lbrace 1,2,...,n(n-1)+1\rbrace$. Let $A_1$ and $A_2$ be the induced subgraphs of $B_n$ such that $A_1=B_n[\lbrace u_1,u_2,...,u_{n-1}\rbrace]$ and $A_2=B_n[\lbrace v_1,v_2,...,v_{n-1}\rbrace]$. \\

			\noindent
			\textbf{\maketitle Construction of an Antimagic Labeling:}\\
			Step 1: We label the edges of the subgraph $A_1$ using $\Big\lbrace 1,2,...,\frac{(n-1)(n-2)}{2}\Big\rbrace$.\\
			
			\noindent
			Step 2: We label the edges of the subgraph $A_2$ using  $\Big\lbrace\frac{(n-1)(n-2)}{2}+1,\frac{(n-1)(n-2)}{2}+2,...,(n-1)(n-2)\Big\rbrace$.  This labeling leads to the partial vertex sums (need not be antimagic) $w'(u_i)$ and $w'(v_i)$, $1\leq i \leq n-1$. The updated $u_i$'s are rewritten as $a_j$'s, $1\leq i,j \leq n-1$ such that $w'(a_j)\leq w'(a_{j+1})$, $1\leq j \leq n-2$ and the updated $v_i$'s are rewritten as $b_j$'s, $1\leq i,j \leq n-1$ such that $w'(b_j)\leq w'(b_{j+1})$, $1\leq j \leq n-2$.  Also, $w'(a_i)<w'(b_j)$, $1 \leq i,j \leq n-1$.\\
			
			\noindent
			Step 3: We label the edges $ua_i$, $vb_j$, $uv$, $1\leq i,j \leq n-1$ using $\lbrace (n-1)(n-2)+1,....(n-1)(n-2)+(n-1)\rbrace$, $\lbrace (n-1)(n-2)+(n-1)+1,...,(n-1)(n-2)+(n-1)+(n-1)\rbrace$ and $n(n-1)+1$ respectively such that,
			\begin{align}
				f(ua_i)=&(n-1)(n-2)+i, 1\leq i \leq n-1\nonumber\\
				f(vb_j)=&(n-1)(n-2)+(n-1)+j, 1\leq j \leq n-1\nonumber\\
				f(uv)=&n(n-1)+1\nonumber	
			\end{align}
			\vspace{-0.5cm}
			We now provide an algorithm to construct an antimagic lableing $f$ of $B_n$\\
			
			\noindent
			\textbf{\maketitle Algorithm:}\\
			
			\noindent
			\textbf{STEP 1:} Label the edges of the subgraphs $A_1$ and $A_2$\\
			1: $f(E(A_1))\leftarrow [1,\frac{(n-1)(n-2)}{2}]$\\
			2: $f(E(A_2))\leftarrow [\frac{(n-1)(n-2)}{2}+1,(n-1)(n-2)]$\\
			
			\noindent
			\textbf{STEP 2:} Label the edges incident with $a_i$, $1\leq i \leq n-1$\\
			3: Sort the partial vertex sums $w'(a_i)$, $1\leq i \leq n-1$ as $w'(a_i)\leq w'(a_{i+1})$, $1\leq i \leq n-2$\nonumber\\
			4: Sort the lables $(n-1)(n-2)+1,(n-1)(n-2)+2,...,(n-1)(n-2)+(n-1)$ in an increasing order\\
			5: \textbf{for} $i=1,2,...,n-1$ \textbf{do}\\
			6: $f(ua_i)\leftarrow (n-1)(n-2)+i$\\
			
			\noindent
			\textbf{STEP 3:} Label the edges incident with $b_j$, $1\leq j \leq n-1$ and an edge $uv$\\
			7: Sort the partial vertex sums $w'(b_j)$, $1\leq j \leq n-1$ as $w'(b_j)\leq w'(b_{j+1})$, $1\leq j \leq n-2$\\
			8: Sort the lables $(n-1)(n-2)+(n-1)+1,(n-1)(n-2)+(n-1)+2,...,(n-1)(n-2)+(n-1)+(n-1)$ in an increasing order\\
			9: \textbf{for} $j=1,2,...,n-1$ \textbf{do}\\
			10: $f(vb_j)\leftarrow (n-1)(n-2)+(n-1)+j$\\
			11: $f(uv)\leftarrow n(n-1)+1$
			\\
			
			\noindent
			\textbf{\maketitle Proof of Antimagicness:}\\
			This labeling leads to the distinctness on the entire vertex sums as follows:
			\begin{align}
				w(a_i)=&w'(a_i)+(n-1)(n-2)+i \nonumber\\
				&< w(a_{i+1})=w'(a_{i+1})+(n-1)(n-2)+(i+1),1\leq i \leq n-2\nonumber\\
				w(b_j)=&w'(b_j)+(n-1)(n-2)+(n-1)+j\nonumber\\  &<w(b_{j+1})=w'(b_{j+1})+(n-1)(n-2)+(n-1)+(j+1),1\leq j \leq n-2\nonumber
			\end{align}
			This clearly shows that $w(a_i)<w(b_j)$, $1\leq i,j \leq n-1$.  
			\begin{align}
				\text{And, }w(u)=&\sum_{i=1}^{n-1}f(ua_i)+n(n-1)+1\nonumber\\
				&< w(v)=\sum_{j=1}^{n-1}f(vb_j)+n(n-1)+1\nonumber
			\end{align} 
			since $\sum_{i=1}^{n-1}f(ua_i) < \sum_{j=1}^{n-1}f(vb_j)$.\\
			Let us define,
			
			\noindent
			set 1: $\big\lbrace \frac{(n-1)(n-2)}{2}+1,....,(n-1)(n-2) \big\rbrace$\\
			set 2: $\lbrace (n-1)(n-2)+(n-1)+1,..., (n-1)(n-2)+(n-1)+(n-1)=n(n-1)\rbrace$\\
			set 3: $\lbrace (n-1)(n-2)+1,...,(n-1)(n-2)+(n-1) \rbrace$\\
			set 4: $n(n-1)+1$\\
			
			The maximum of $w(b_j)$, $1\leq j \leq n-1$ is $w(b_{n-1})$ which is the sum of any $n-2$ labels of set 1 and any one label of set 2; and $w(u)$ is the sum of all labels of set 3 and set 4. Observe that $d(b_{j})<d(u)$, $1\leq j \leq n-1$. Clearly, $w(b_j) < w(u)$, $1\leq j \leq n-1$. Hence, all the vertices are distinct since $w(a_i) < w(b_j) < w(u) < w(v)$, $1 \leq i,j \leq n-1$.\\
			
			\noindent
			\textbf{\maketitle Time Complexity:}\\
			The assignments in Step 1 takes constant time. Since, the partial vertex sums are sorted in line 3 of Step 2, it requires $O(nlogn)$ time. Again, Step 3 requires the time $O(nlogn)$  due to the same fact that the partial vertex sums are sorted in line 7. Hence, the time complexity of the above
			algorithm is $O(nlogn)$.
		\end{proof}
	\end{theorem}
	\noindent
	Next, we illustrate the above labeling process in Figure \ref{fig:2} for the $4$-barbell graph, $B_4$ .
	
	Note that $w(a_1)=10$, $w(a_2)=12$, $w(a_3)=14$, $w(u)=37$ and $w(b_1)=19$, $w(b_2)=21$, $w(b_3)=23$, $w(v)=46$.
	\subsection{Edge Corona of $B_{x,n}$ and $H$}
\begin{figure}[ht!]
		\centering
		\includegraphics[width=0.9\linewidth, height=0.35\textheight]{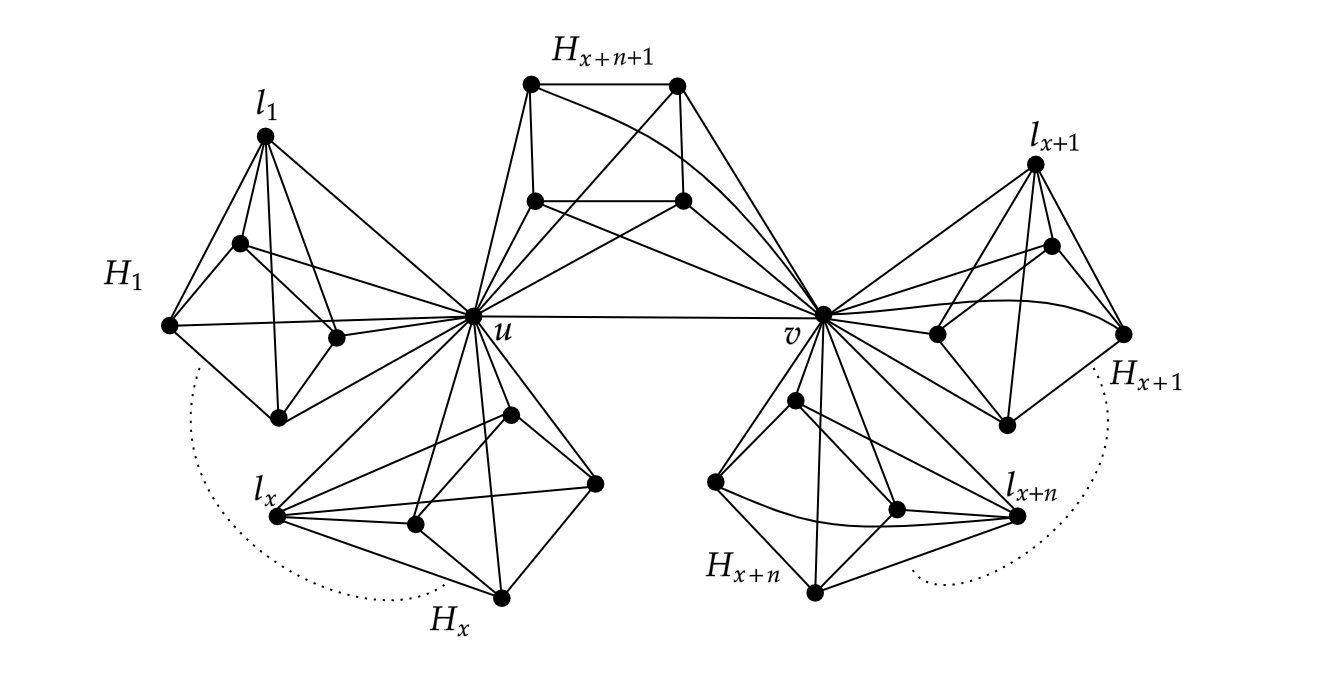}
		\caption{A representation of $B_{x,n}\diamond C_4$.}
		\label{fig:3}
	\end{figure}

	Let the vertices $\lbrace l_1,l_2,...,l_x\rbrace$ and $\lbrace l_{x+1},l_{x+2},...,l_{x+n}\rbrace$ be adjacent to the apex vertices $u$ and $v$ of $K_{1,x}$ and $K_{1,n}$  respectively to form a bistar graph $B_{x,n}$ with $uv\in E(B_{x,n})$ where $x,n\geq 2$. Note that the graph $B_{x,n}$ contains $x+n+2$ vertices and $x+n+1$ edges. To construct the graph $B_{x,n}\diamond H$, we require one copy of $B_{x,n}$ and $|E(B_{x,n})|$ copies of $H$ namely $H_1,H_2,...,H_{x+n+1}$ (each $H_i$ is a $k$-regular graph on $m$ vertices). Let the  edge $l_iu$ be complete to $H_i$, $1\leq i \leq x$ and the edge $l_iv$ be complete to $H_i$, $x+1\leq i \leq x+n$. And, the edge $uv$ is complete to $H_{x+n+1}$.
	The graph $\bigcup\limits_{i=1}^{x+n+1} H_i$ contains $(x+n+1)m$ vertices and $(x+n+1)\frac{mk}{2}$ edges. Therefore, the graph $B_{x,n}$ $\diamond$ $H$ contains $x+n+2+m(x+n+1)$ vertices and  $(x+n+1)+2m(x+n+1)+\frac{mk}{2}(x+n+1)=z$ (say) edges. For a better understanding of the graph $B_{x,n}\diamond H$, see Figure \ref{fig:3} ($H$ isomorphic to $C_4$).
	The left dotted curve represent the graph $\bigcup\limits_{i=2}^{x-1}(H_i+\lbrace l_i\rbrace)$. Similarly, the right dotted curve represent the graph $\bigcup\limits_{i=x+2}^{x+n-1}(H_i+\lbrace l_i\rbrace)$.
	
\begin{figure}[ht!]
		\centering
		\includegraphics[width=1\linewidth, height=0.35\textheight]{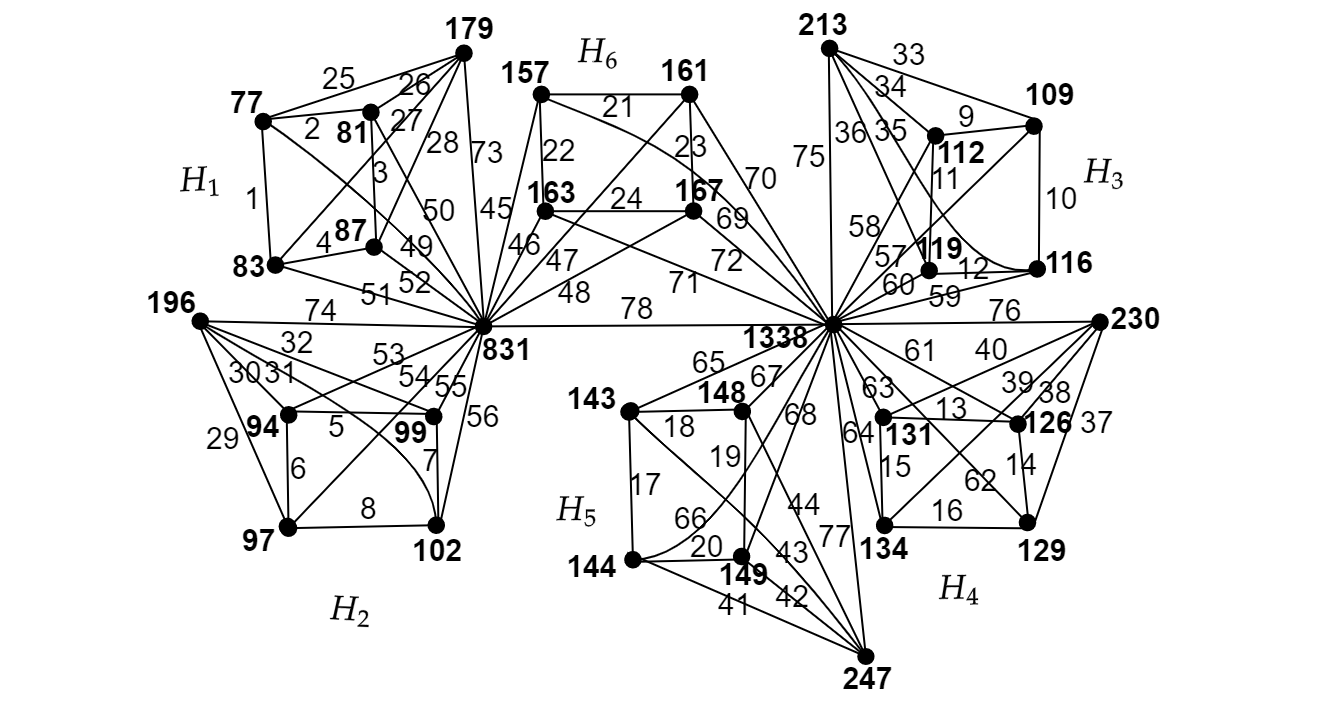}
		\caption{An antimagic labeling of $B_{2,3}\diamond C_4$.}
		\label{fig:4}
	\end{figure}

	\begin{theorem}\label{Theorem 3.2}
		$B_{x,n}$ $\diamond$ $H$ is antimagic for $ x,n\geq 2$, $x\leq n$ and $H$ is a connected $k$-regular graph on $m\geq 2$ vertices.
		\begin{proof}
			Let $\lbrace u_1,u_2,...,u_{mx}\rbrace$ be the vertices of $\bigcup \limits_{i=1}^{x}H_i$ that are adjacent to $u$ and $\lbrace v_1,v_2,...,v_{mn}\rbrace$ be the vertices of $\bigcup \limits_{i=x+1}^{x+n}H_i$ that are adjacent to $v$.
			Let $\lbrace w_1,w_2,...,w_m\rbrace$ be the vertices of $H_{x+n+1}$ that are adjacent to  both $u$ and $v$. \\
			
			\noindent
			\textbf{\maketitle{Construction of an Antimagic Labeling:}}\\
			\indent	Step 1: We label the edges of graphs $H_1,H_2,...,H_x$ in an order from smallest to the largest label available in the set $\lbrace1,2,...,x(\frac{mk}{2})\rbrace$.\\
			
			Step 2: We label the edges of graphs $H_{x+1},H_{x+2},...,H_{x+n}$ in an order from smallest to the largest label available in the set $\lbrace x(\frac{mk}{2})+1,x(\frac{mk}{2})+2,...,x(\frac{mk}{2})+n(\frac{mk}{2})\rbrace$.\\
			
			Step 3: We label the edges of graph $H_{x+n+1}$ using $\Big\lbrace x(\frac{mk}{2})+n(\frac{mk}{2})+1,x(\frac{mk}{2})+n(\frac{mk}{2})+2,...,x(\frac{mk}{2})+n(\frac{mk}{2})+\frac{mk}{2}=g\Big\rbrace$. \\
			
			Step 4: We label the set of edges incident with $l_1,l_2,...,l_x$ excluding the edges  $ul_1,ul_2,...,ul_x$ in an order from smallest to the largest label available in the set $\lbrace g+1,g+2,...,g+xm\rbrace$ respectively . \\ 
			
			Step 5: We label the set of edges incident with $l_{x+1},l_{x+2},...,l_{x+n}$ excluding the edges $vl_{x+1},vl_{x+2},...,vl_{x+n}$ in an order from smallest to the largest label available in the set $\lbrace g+xm+1,g+xm+2,...,g+xm+nm\rbrace$ respectively.\\
			
			Step 6: We label the edges $uw_1,uw_2,...,uw_m$ using  $\lbrace g+xm+nm+1,g+xm+nm+2,...,g+xm+nm+m=h\rbrace$.  With respect to the above labeling we obtain the partial vertex sums (need not be antimagic) $w'(u_i)$, $ 1\leq i \leq mx$, $w'(v_i)$, $ 1\leq i \leq mn$, $w'(w_i)$, $1\leq i \leq m$ and $w'(l_i)$, $1\leq i \leq x+n$.\\
			
			We exclude the partial vertex sum of the vertex $u$. Merging all the above partial vertex sums, we update the new vertex sums as  $w'(a_i)\leq w'(a_{i+1})$, $1\leq i \leq  x(m+1)+n(m+1)+m-1$.\\ 
			
			Step 7: We label the edge $uv$ as $f(uv)=z$.\\
			
			Step 8: Finally, we label the remaining edges using $\lbrace h+1,h+2,...,h+x(m+1)+n(m+1)+m\rbrace$ in such a way that 
			$f(sa_i)=h+i$, $1\leq i \leq x(m+1)+n(m+1)+m$ where $s=u$ for all $a_i$ that are adjacent to $u$ and $s=v$ for all $a_i$ that are adjacent to $v$.
			
			We now provide an algorithm to construct an antimagic labeling $f$ of $B_{x,n}\diamond H$\\
			
			\noindent
			\textbf{\maketitle Algorithm:}\\
			
			\noindent
			\textbf{STEP 1:} Label the edges of the graphs $H_i$, $i=1,2,...,x,x+1,...,x+n,x+n+1$\\
			1: \textbf{for} $i=1,2,...,x,x+1,...,x+n,x+n+1$ \textbf{do}\\
			2: $f(E(H_i))\leftarrow\big[(i-1)\frac{mk}{2}+1,i\frac{mk}{2}\big]$\\
			
			\noindent
			\textbf{STEP 2:} Label the edges incident with $l_1,l_2,...,l_{x},l_{x+1},...,l_{x+n}$ excluding the edges $ul_1,ul_2,...,ul_x,ul_{x+1},...,ul_{x+n}$ respectively\\
			3: \textbf{for} $i=1,2,...,x,x+1,...,x+n$ \textbf{do}\\
			4: $f(E(l_i))\leftarrow [g+(i-1)m+1,g+im]$ excluding $f(ul_i)$\\
			
			\noindent
			\textbf{STEP 3:} Label the edges $uw_i$, $1\leq i \leq m$ and $uv$\\
			5: \textbf{for} $i=1,2,...,m$ \textbf{do}\\
			6: $f(uw_i)\leftarrow g+xm+nm+i$\\
			7: $f(uv)\leftarrow z$\\
			
			\noindent
			\textbf{STEP 4:} Label the remaining edges of $B_{x,n}\diamond H$\\
			8: Sort the partial vertex sums $w'(a_i)$, $1\leq i\leq x(m+1)+n(m+1)+m$ as $w'(a_i)\leq w'(a_{i+1})$, $1\leq
			i\leq x(m+1)+n(m+1)+m-1$\\
			9: Sort the available labels $h+1,h+2,...,h+x(m+1)+n(m+1)+m$ in an increasing order\\
			10: \textbf{for} $i=1,2,...,x(m+1)+n(m+1)+m$ \textbf{do}\\
			11: $f(sa_i)\leftarrow h+i$, $s=u$ $\forall$ $a_i$ adjacent to $u$, $s=v$ $\forall$ $a_i$ adjacent to $v$\\
			
			\noindent
			\textbf{\maketitle Proof of Antimagicness:}\\
			This labeling leads to the distinctness on the entire vertex sums as follows:
			\begin{align}
				w(a_i)&=w'(a_i)+h+i <
				w(a_{i+1})=w'(a_{i+1})+h+(i+1)\nonumber\\ 
				&\text{for }1\leq i \leq x(m+1)+n(m+1)+m-1\nonumber 
			\end{align}
			Let us define, \\
			\noindent
			set 1: $\lbrace g+xm+nm+1,...,g+xm+nm+m=h\rbrace, \lbrace h+1,h+2,...,h+xm\rbrace, \lbrace h+xm+nm+m+1,h+xm+nm+2,...,h+xm+nm+x\rbrace, \lbrace z\rbrace$\\
			\noindent\\
			set 2: $\lbrace h+xm+1,h+xm+2,...,h+xm+nm\rbrace, \lbrace h+xm+nm+1,h+xm+nm+2,...,h+xm+nm+m\rbrace,\lbrace h+xm+nm+m+x+1,h+xm+nm+m+x+2,...,h+xm+nm+m+x+n\rbrace, \lbrace z\rbrace$ \\
			\indent The edges incident with $u$ receives the labels of set 1 and the edges incident with $v$ receives the labels of set 2. Observe that $d(u)\leq d(v)$. From the above, it is clear that the sum of all the labels of set 1 is less than the sum of all the labels of set 2. So, it is clearly shown that $w(u)<w(v)$ using the labels of set 1 and set 2. 
			And the maximum of $w(a_i)$ is $w(a_{x(m+1)+n(m+1)+m})$ which is the sum of $w'(a_{x(m+1)+n(m+1)+m})$ and $z-1$. Also,
			\begin{align}
				w(a_{x(m+1)+n(m+1)+m})&=w'(a_{x(m+1)+n(m+1)+m})+(z-1) <\nonumber\\ w(u)&=w'(u)+z+\sum\limits_{i=1}^{xm}(h+i)+\sum\limits_{i=1}^x(h+xm+nm+m+i)\nonumber\\ 
				\text{where } w'(a_{x(m+1)+n(m+1)+m})&=\sum\limits_{i=1}^{m}(g+xm+(n-1)m+i) \nonumber\\ 
				\text{and }w'(u)&=\sum\limits_{i=1}^m(g+xm+nm+i)\nonumber
			\end{align}
			Hence, $w(a_i)<w(u)<w(v)$, $1\leq i \leq x(m+1)+n(m+1)+m$.\\
			
			\noindent
			\textbf{\maketitle Time Complexity:}\\
			The assignments in Step 1 reaches at most $x+n+1$ times to label the edges of the graphs $H_i$ and so, it takes $O(x+n)$ time. The same $O(x+n)$ time is required for Step 2 to label the edges incident with the vertices $l_i$ excluding the edges $ul_i$ respectively. An assignment in line 6 of Step 3 reaches at most $m$ vertices and hence it require $O(m)$ time. Line 8 in Step 4 requires $O(ylogy)$ time since the partial vertex sums are sorted where $y=x(m+1)+n(m+1)+m$. Hence, the time complexity of the above algorithm is $O(ylogy)$.\end{proof}
	\end{theorem}
	Next, we illustrate the above labeling process in Figure \ref{fig:4} for the graph $B_{2,3}\diamond C_4$. Note that $w(u)=831$, $w(v)=1338$ and $w(a_1)=77$, $w(a_2)=81$, $w(a_3)=83$, $w(a_4)=87$, $w(a_5)=94$, $w(a_6)=97$, $w(a_7)=99$, $w(a_8)=102$, $w(a_9)=109$, $w(a_{10})=112$, $w(a_{11})=116$, $w(a_{12})=119$, $w(a_{13})=126$, $w(a_{14})=129$, $w(a_{15})=131$, $w(a_{16})=134$, $w(a_{17})=143$, $w(a_{18})=144$, $w(a_{19})=148$, $w(a_{20})=149$, $w(a_{21})=157$, $w(a_{22})=161$, $w(a_{23})=163$, $w(a_{24})=167$, $w(a_{25})=179$, $w(a_{26})=196$, $w(a_{27})=213$, $w(a_{28})=230$, $w(a_{29})=247$.

	\begin{corollary}\rm\label{Corollary 1}
		$B_{x,n}$ $\diamond$ $H$ admits an antimagic labeling for $ x,n\geq 2$ and $x>n$ where $H$ is a connected $k$-regular graph on $m\geq 2$ vertices, $k\geq 1$.
		\begin{proof}
			In Theorem \ref{Theorem 3.2}, we proved that $B_{x,n}$ $\diamond$ $H$, $x\leq n$ is antimagic. Obviously $B_{x,n}$ $\diamond$ $H$, $x\geq n$ is isomorphic to $B_{x,n}$  $\diamond$ $H$, $x\leq n$ and hence the result.
		\end{proof}
	\end{corollary}
	\subsection{Edge Corona of $C_m$ and $C_n$}
	To construct the graph $C_m\diamond C_n$, $m,n\geq3$ we join one copy of $C_m$ with $|E(C_m)|$ copies of $C_n$ namely $C_{n_{1}},C_{n_{2}},...,C_{n_{m}}$ such that the end vertices of $i^{th}$ edge of $C_m$ is adjacent to every vertex in the $i^{th}$ copy of $C_n$. The following lemma is already proved in \cite{13}. We again prove the lemma with different labels so as to use it in Theorem \ref{Theorem 3.3}
	\begin{lemma}\rm\label{3.6}
		Cycle $C_m$, $m\geq 3$ is antimagic.
		\begin{proof}
			Let $V(C_m)=\lbrace v_1,v_2,...,v_m\rbrace$ and $E(C_m)=\lbrace v_1v_2, v_1v_3 \rbrace\bigcup\lbrace v_{i}v_{i+2}\vspace{0.1cm}\big|\vspace{0.1cm} i=2,...,m-2\rbrace\bigcup\lbrace v_{m-1}v_m\rbrace$. Now, we shall show that  $C_m$ is antimagic with different set of labels. In the proof of the lemma in \cite{13}, adding $3mn$ to the edge labels and $6mn$ to the vertex sums, we get $f(v_1v_2)=3mn+1$, $f(v_1v_3)=3mn+2$, $f(v_iv_{i+2})=3mn+i+1$, for $2\leq i \leq m-2$, and $f(v_{m-1}v_m)=3mn+m$. The edge labeling induces the following ordering on vertices as $w(v_1)<w(v_2)<...<w(v_m)$ since the vertex sums are 
			
			$$w(v_i)= 
			\begin{cases}
				6mn+3&\text{ if }i =1;\\
				6mn+2i&\text{ if }i=2,...,m-1;\\
				6mn+2m-1&\text{ if }i=m\\
			\end{cases}$$
			Hence, $C_m$ is antimagic.
		\end{proof}
	\end{lemma}
	\begin{lemma}\rm\label{3.7}
		Cycles $C_{n_{j}}$, $1\leq j\leq m$ are antimagic.
		\begin{proof}
			Let $V(C_{n_{j}})=\lbrace u^j_1,u^j_2,...,u^j_n \rbrace$,  $1\leq j \leq m$ and $E(C_{n_{j}})= \lbrace u^j_1u^j_2,u^j_1u^j_3\rbrace\\
			\bigcup\lbrace u^j_iu^j_{i+2}\hspace{0.1cm}|\hspace{0.1cm} u_i=2,...,n-2\rbrace\bigcup\lbrace u^j_{n-1}u^j_n\rbrace$, $1\leq j \leq m$. Now, we shall show that $C_{n_{j}}$ is antimagic with different set of labels. In the proof of the lemma in \cite{13}, adding $(j-1)n$ to the edge labels and $2(j-1)n$ to the vertex sums, we get $f(u_1^ju_2^j)=(j-1)n+1$, $f(u_1^ju_3^j)=(j-1)n+2$, $f(u_i^ju_{i+2}^j)=(j-1)n+i+1$, for $2\leq i \leq n-2$, and
			$f(u^j_{n-1}u^j_n)=(j-1)n+n$. Note that, $w(u_1^j)<w(u_2^j)<...<w(u_m^j)$ since the vertex sums are 
			$$w(u_i^j)= 
			\begin{cases}
				2(j-1)n+3&\text{ if }i=1;\\
				2(j-1)n+2i&\text{ if }i=2,...,n-1;\\
				2(j-1)n+2n-1&\text{ if }i=n\\
			\end{cases}$$
			Hence, $C_{n_{j}}$, $1\leq j\leq m$ are antimagic
		\end{proof}
	\end{lemma}
	
\begin{figure}[ht!]
		\centering
		\includegraphics[width=0.8\linewidth, height=0.35\textheight]{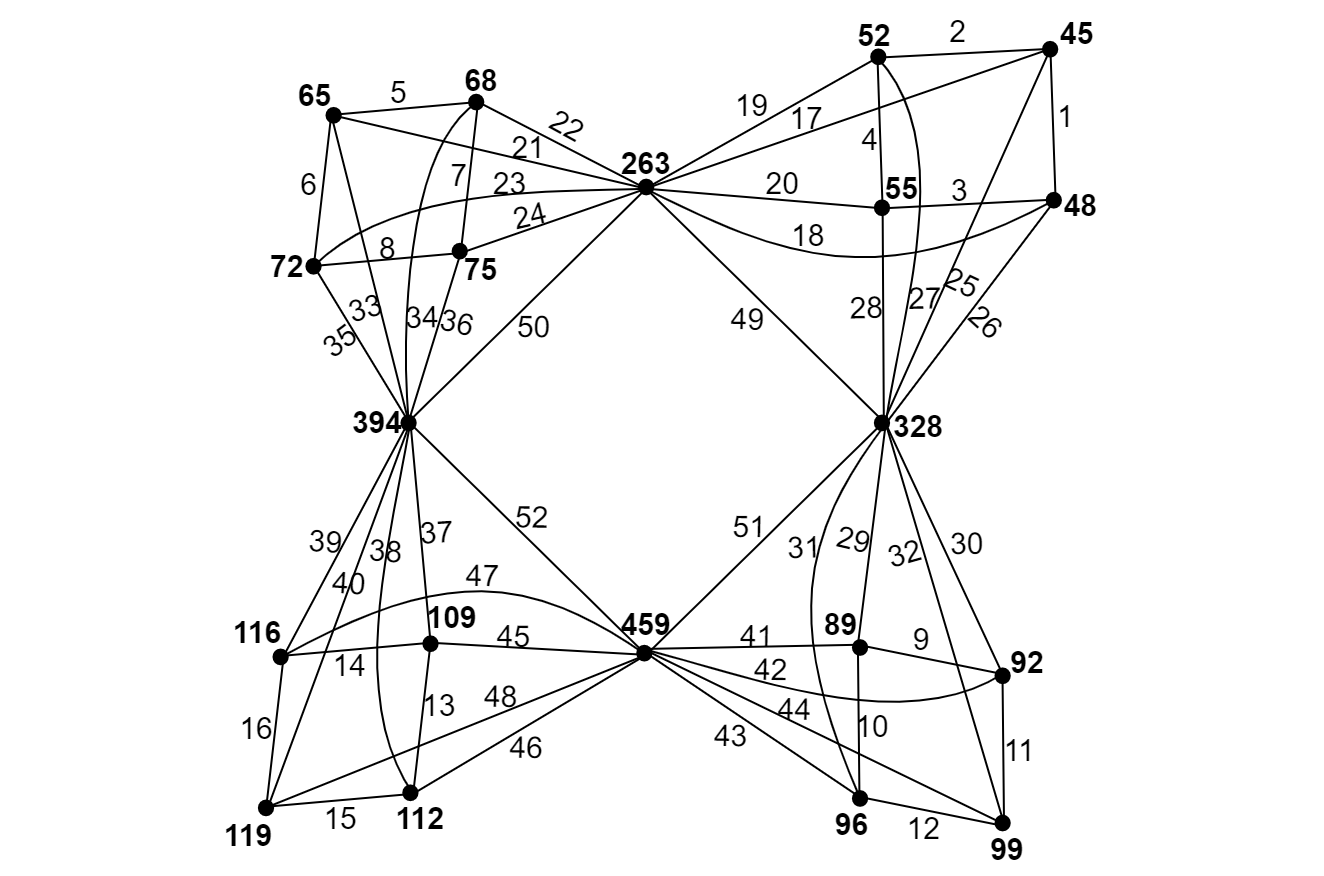}
		\caption{An antimagic labeling of $C_4 \diamond C_4$.}
		\label{fig:5}
	\end{figure}

	\begin{theorem}\label{Theorem 3.3}
		$C_m\diamond C_n$ is antimagic for $m,n\geq3$ 
		\begin{proof}
			Let the vertex set, edge set, and an antimagic labeling for the graphs $C_m$ and $m$ copies of $C_n$ be defined as before. Rename the vertex sums $w(u_i^j)$, $w(v_i)$ as $w'(u_i^j)$ and $w'(v_i)$ respectively (the vertex sums of $C_m$ and $C_{n_j}$ are considered as partial vertex sums of $C_m \diamond C_n$). The adjaceny in  $C_m \diamond C_n$ is defined as follows:\\						  \textbf{Case i.} when $m$ is even, $m=2k$, $k\in 
			\mathbb{Z}^+-\lbrace1\rbrace$\\
			$v_1 v_2$ is complete to $C_{n_1}$\\
			$ v_1v_3\text{ is complete to } C_{n_2}$\\
			$v_2v_4\text{ is complete to }C_{n_3}$\\
			$v_3v_5\text{ is complete to }C_{n_4}$ \\
			$\vdots$\\
			$v_{2k-2}v_{2k}\text{ is complete to }C_{n_{2k-1}}$\\
			$ v_{2k}v_{2k-1}\text{ is complete to }C_{n_{2k}}$ \\	
			
			\noindent   \textbf{Case ii.}   when $m$ is odd, $m=2k+1$, $k\in \mathbb{Z}^+$\\
			$v_1v_2\text{ is complete to }C_{n_1} $ \\
			$v_1v_3\text{ is complete to }C_{n_2} $\\
			$ v_2v_4\text{ is complete to }C_{n_3}$\\
			$v_3v_5\text{ is complete to }C_{n_4}$\\
			$\vdots$\\
			$v_{2k-1}v_{2k+1}\text{ is complete to }C_{n_{2k}}$\\
			$v_{2k+1}v_{2k}\text{ is complete to }C_{n_{2k+1}}	$\\	
			
			\noindent 
			\textbf{\maketitle Construction of an Antimagic Labeling:}\\
			The edge labels of $C_m\diamond C_n$ are $\lbrace1,2,...,3mn+m\rbrace$. As in the above lemma \ref{3.6},\ref{3.7} we assign the labels  $\lbrace1,2,...,mn\rbrace$ to the edges of $C_{n_j}$, $1\leq j\leq m$ and the labels $\lbrace 3mn+1,...,3mn+m\rbrace$  to the edges of $C_m$. Note that the induced edges between $C_m$ and $C_{n_{j}}$, $1\leq j \leq m$ are to be labelled with  $\lbrace mn+1,mn+2,...,3mn\rbrace$ such that,
			\noindent
			\begin{align}
				f(v_1u^1_i)&=mn+i, 1\leq i \leq n\nonumber\\
				f(v_1u^2_i)&=mn+n+i, 1\leq i \leq n\nonumber\\
				f(v_2u^1_i)&=mn+2n+i, 1\leq i \leq n\nonumber\\
				f(v_2u^3_i)&=mn+3n+i, 1\leq i \leq n\nonumber\\
				f(v_3u^2_i)&=mn+4n+i, 1\leq i \leq n\nonumber\\
				f(v_3u^4_i)&=mn+5n+i, 1\leq i \leq n\nonumber\\
				f(v_4u^3_i)&=mn+6n+i, 1\leq i \leq n\nonumber\\
				f(v_4u^5_i)&=mn+7n+i, 1\leq i \leq n\nonumber\\
				f(v_5u^4_i)&=mn+8n+i, 1\leq i \leq n\nonumber\\
				f(v_5u^6_i)&=mn+9n+i, 1\leq i \leq n\nonumber\\
				&\vdots\nonumber\\
				f(v_mu^{m-1}_i)&=mn+(m-1)2n+i, 1\leq i \leq n\nonumber\\
				f(v_mu^m_i)&=3mn-n+i, 1\leq i \leq n\nonumber
			\end{align}
			
			\noindent
			\textbf{\maketitle Algorithm:}\\
			
			\noindent
			\textbf{STEP 1:} Label the edges of $C_m$\\
			1: $f(v_1v_2)\leftarrow3mn+1$\\
			2: $f(v_1v_3)\leftarrow3mn+2$\\
			3: \textbf{for} $i=2$ to $m-2$ \textbf{do}\\
			4: $f(v_iv_{i+2})\leftarrow 3mn+i+1$\\
			5: $f(v_{m-1}v_m)\leftarrow 3mn+m$\\
			
			\noindent
			\textbf{STEP 2:} Label the edges of $C_{n_j}$, $1\leq j \leq m$\\
			6: \textbf{for} $j=1,2,...,m$, $i=2$ to $n-2$   \textbf{do}\\
			7: $f(u_1^ju_2^j)\leftarrow (j-1)n+1$\\
			8: $f(u_1^ju_3^j)\leftarrow (j-1)n+2$\\
			9: $f(u_i^ju_{i+2}^j)\leftarrow (j-1)n+i+1$\\
			10: $f(u_{n-1}^ju_n^j)\leftarrow (j-1)n+n$\\
			
			\noindent
			\textbf{STEP 3:} Label the induced edges between $C_{n_j}$, $1\leq j \leq m$ and $C_m$\\
			11: \textbf{for} $i=1$ to $n$ \textbf{do}\\
			12: $f(v_1u^1_i)\leftarrow mn+i$\\
			$f(v_1u^2_i)\leftarrow mn+n+i$\\
			$f(v_2u^1_i)\leftarrow mn+2n+i$\\
			$f(v_2u^3_i)\leftarrow mn+3n+i$\\
			$f(v_3u^2_i)\leftarrow mn+4n+i$\\
			$f(v_3u^4_i)\leftarrow mn+5n+i$\\
			$f(v_4u^3_i)\leftarrow mn+6n+i$\\
			$f(v_4u^5_i)\leftarrow mn+7n+i$\\
			$f(v_5u^4_i)\leftarrow mn+8n+i$\\
			$f(v_5u^6_i)\leftarrow mn+9n+i$\\
			\vdots\nonumber\\
			$f(v_mu^{m-1}_i)\leftarrow mn+(m-1)2n+i, 1\leq i \leq n$\\
			$f(v_mu^m_i)\leftarrow3mn-n+i, 1\leq i \leq n$\\
			
			\noindent
			\textbf{\maketitle Proof of Antimagicness:}
			This labeling leads to the distinctness on the entire vertex sums as follows:
			\noindent
			\begin{align}
				&w(u^1_{i_1})=w'(u^1_{i_1})+f(v_1u^1_{i_1})+f(v_2u^1_{i_1})\nonumber\\
				&<	w(u^2_{i_2})=w'(u^2_{i_2})+f(v_1u^2_{i_2})+f(v_3u^2_{i_2})\nonumber\\
				&<	w(u^3_{i_3})=w'(u^3_{i_3})+f(v_2u^3_{i_3})+f(v_4u^3_{i_3})\nonumber\\
				&<	w(u^4_{i_4})=w'(u^4_{i_4})+f(v_3u^4_{i_4})+f(v_5u^4_{i_4})\nonumber\\
				&\vdots\nonumber\\
				&<	w(u^{m-1}_{i_{m-1}})=w'(u^{m-1}_{i_{m-1}})+f(v_{m-2}u^{m-1}_{i_{m-1}})+f(v_mu^{m-1}_{i_{m-1}})\nonumber\\
				&<	w(u^m_{i_m})=w'(u^m_{i_m})+f(v_mu^m_{i_m})+f(v_{m-1}u^m_{i_m}), 1\leq i_1,i_2,...,,i_m \leq n
			\end{align}
			And, 
			\begin{align}
				&w(u^1_{i_1})=w'(u^1_{i_1})+f(v_1u^1_{i_1})+f(v_2u^1_{i_1})\nonumber\\
				&<w(u^1_{i_1+1})=w'(u^1_{i_1+1})+f(v_1u^1_{i_1+1})+f(v_2u^1_{i_1+1});\nonumber\\
				&w(u^2_{i_2})=w'(u^2_{i_2})+f(v_1u^2_{i_2})+f(v_3u^2_{i_2})\nonumber\\
				&<w(u^2_{i_2+1})=w'(u^2_{i_2+1})+f(v_1u^2_{i_2+1})+f(v_3u^2_{i_2+1});\nonumber\\
				&w(u^3_{i_3})=w'(u^3_{i_3})+f(v_2u^3_{i_3})+f(v_4u^3_{i_3})\nonumber\\
				&<w(u^3_{i_3+1})=w'(u^3_{i_3+1})+f(v_2u^3_{i_3+1})+f(v_4u^3_{i_3+1})\nonumber\\
				&\vdots\nonumber\\
				&w(u^m_{i_m})=w'(u^m_{i_m})+f(v_mu^m_{i_m})+f(v_{m-1}u^m_{i_m})\nonumber\\
				&<w(u^m_{i_m+1})=w'(u^m_{i_m+1})+f(v_mu^m_{i_m+1})+f(v_{m-1}u^m_{i_m+1}), 1\leq i_1,i_2,...,,i_m \leq n-1\nonumber	
			\end{align}\\
			\noindent
			Also,
			\begin{align}
				w(v_i)=&w'(v_i)+\sum\limits_{j=1}^{2n}(mn+(i-1)2n+j)\nonumber\\&<w(v_{i+1})=w'(v_{i+1})+\sum\limits_{j=1}^{2n}(mn+(i)2n+j), 1\leq i \leq m-1\nonumber
			\end{align}
			Let us define, \\
			\noindent
			set 1: $\lbrace 1,2,...,mn \rbrace$\\
			set 2: $\lbrace mn+1,...,3mn\rbrace$\\
			set 3: $\lbrace 3mn+1,...,3mn+m\rbrace$\\
			
			Observe that $w(u^i_j)<w(v_s)$, $1\leq i,s \leq m$, $1\leq j \leq n$. Since, $w(u^i_j)$ is the sum of any two labels of the set 1 and any two labels of the set 2; $w(v_s)$ is the sum of any two labels of set 3 and any $2n$ labels of set 2, $w(u^i_j)<w(v_s)$, $1\leq i,s \leq m$, $1\leq j \leq n$. Hence the vertex sums of the graph $C_m \diamond C_n$ are distinct.\\
			
			\noindent
			\textbf{\maketitle Time Complexity:}\\
			Step 1 takes $O(m)$ time since the line 4 iterates at most $m-3$ times. Step 2 requires $O(mn)$ time since the line 9 iterates for $m(n-3)$ times. Finally, the step 3 takes $O(mn)$ time since it iterates $2mn$ times. Hence, the time complexity of the above algorithm is $O(mn)$. 
		\end{proof}
	\end{theorem}
	\noindent
	We illustrate the above labeling process in Figure \ref{fig:5} for the graph $C_4 \diamond C_4$. Note that $w(u^1_1)=45$, $w(u^1_2)=48$, $w(u^1_3)=52$, $w(u^1_4)=55$, $w(u^2_1)=65$, $w(u^2_2)=68$, $w(u^2_3)=72$, $w(u^2_4)=75$, $w(u^3_1)=89$, $w(u^3_2)=92$, $w(u^3_3)=96$, $w(u^3_4)=99$, $w(u^4_1)=109$, $w(u^4_2)=112$, $w(u^4_3)=116$, $w(u^4_4)=119$ and $w(v_1)=263$, $w(v_2)=328$, $w(v_3)=394$, $w(v_4)=459$.
	
	\section{Conclusions}
	As the conjecture due to Hartsfield and Ringel remains open for all these years we have investigated the antimagic labeling of the barbell graph and edge corona of some classes of graphs. It is also interesting to work on antimagic labeling of the generalized edge corona of graphs.

\end{document}